\documentclass[12pt,a4paper]{article}
\usepackage{amsmath}
\usepackage{enumerate}
\usepackage{amsthm}
\usepackage{theoremref}
\usepackage{amssymb}
\usepackage{t1enc}
\usepackage[latin1]{inputenc}
\usepackage[english]{babel}
\usepackage{amsfonts}
\usepackage{latexsym}
\usepackage{bm}
\usepackage{lineno}
\newtheorem{theorem-section}{Theorem}[section]
\newtheorem{proposition-section}[theorem-section]{Proposition}
\newtheorem{lemma-section}[theorem-section]{Lemma}
\newtheorem{corollary-section}[theorem-section]{Corollary}
\newtheorem{conjecture-section}[theorem-section]{Conjecture}
\newtheorem{claim-section}[theorem-section]{Claim}
\newtheorem{fact-section}[theorem-section]{Fact}
\newtheorem{remark-section}[theorem-section]{Remark}
\newtheorem{observation-section}[theorem-section]{Observation}
\newtheorem{definition-section}[theorem-section]{Definition}
\newtheorem{construction-section}[theorem-section]{Construction}
\newtheorem{problem-section}[theorem-section]{Problem}
\newtheorem{example-section}[theorem-section]{Example}

\newtheorem{theorem-subsection}{Theorem}[subsection]
\newtheorem{proposition-subsection}[theorem-subsection]{Proposition}
\newtheorem{lemma-subsection}[theorem-subsection]{Lemma}
\newtheorem{corollary-subsection}[theorem-subsection]{Corollary}
\newtheorem{conjecture-subsection}[theorem-subsection]{Conjecture}
\newtheorem{claim-subsection}[theorem-subsection]{Claim}
\newtheorem{fact-subsection}[theorem-subsection]{Fact}
\newtheorem{remark-subsection}[theorem-subsection]{Remark}
\newtheorem{observation-subsection}[theorem-subsection]{Observation}
\newtheorem{definition-subsection}[theorem-subsection]{Definition}
\newtheorem{construction-subsection}[theorem-subsection]{Construction}
\newtheorem{problem-subsection}[theorem-subsection]{Problem}
\newtheorem{example-subsection}[theorem-subsection]{Example}

\begin{document}


\title{Irregular independence and irregular domination}

\author{Peter Borg\\[5mm]
Department of Mathematics \\
University of Malta\\
Malta\\
\texttt{peter.borg@um.edu.mt} 
\and 
Yair Caro \\ [5mm]
Department of Mathematics \\
University of Haifa-Oranim \\
Israel \\
\texttt{yacaro@kvgeva.org.il}
\and \\
Kurt Fenech \\ [5mm]
Department of Mathematics \\
University of Malta\\
Malta\\
\texttt{kurt.fenech.10@um.edu.mt}
}

\date{}

\maketitle

\begin{abstract}
If $A$ is an independent set of a graph $G$ such that the vertices in $A$ have different degrees, then we call $A$ an \emph{irregular independent set of} $G$. If $D$ is a dominating set of $G$ such that the vertices that are not in $D$ have different numbers of neighbours in $D$, then we call $D$ an \emph{irregular dominating set of} $G$. The size of a largest irregular independent set of $G$ and the size of a smallest irregular dominating set of $G$ are denoted by $\alpha_{ir}(G)$ and $\gamma_{ir}(G)$, respectively. We initiate the investigation of these two graph parameters. For each of them, we obtain sharp bounds in terms of basic graph parameters such as the order, the size, the minimum degree and the maximum degree, and we obtain Nordhaus--Gaddum-type bounds. We also establish sharp bounds relating the two parameters. Furthermore, we characterize the graphs $G$ with $\alpha_{ir}(G)=1$, we determine those that are planar, and we determine those that are outerplanar.
\end{abstract}

\section{Introduction}
\pagenumbering{arabic}

In this paper, we will consider the notions of irregular independence and irregular domination as counterparts of the notions of regular independence and regular domination (also referred to as fair domination), which were recently introduced in \cite{CHH,CHP}. The formal definitions of these two parameters are as follows.

If $A$ is an independent set of a graph $G$ such that the vertices in $A$ have different degrees, then we call $A$ an \emph{irregular independent set of} $G$. The size of a largest irregular independent set of $G$ will be called the \emph{irregular independence number of $G$} and will be denoted by $\alpha_{ir}(G)$. If $A$ is an independent set of a graph $G$ such that the vertices in $A$ have the same degree, then $A$ is called a \emph{regular independent set of} $G$. The size of a largest regular independent set of $G$ is called the \emph{regular independence number of $G$} and is denoted by $\alpha_{reg}(G)$.  

For a vertex $v$ of a graph $G$, let $N(v)$ denote the set of neighbours of $v$. If $D$ is a dominating set of $G$ such that $|N(u) \cap D| \neq |N(v) \cap D|$ for every two distinct vertices $u$ and $v$ in $V(G) \backslash D$, then we call $D$ an \emph{irregular dominating set of} $G$. The size of a smallest irregular dominating set of $G$ will be called the \emph{irregular domination number of $G$} and will be denoted by $\gamma_{ir}(G)$. If $D$ is a dominating set of $G$ such that $|N (u) \cap D| = |N(v) \cap D|$ for every two vertices $u$ and $v$ in $V(G) \backslash D$, then $D$ is called a \emph{regular dominating set of} $G$. The size of a smallest regular dominating set of $G$ is called the \emph{regular domination number of $G$} and is denoted by $\gamma _ {reg}(G)$. Observe that the notion of irregular domination is an extreme case of the well-studied notion of location-domination \cite{BFH}: a set $D$ is called a \emph{locating-dominating set of} $G$ if $D$ is a dominating set of $G$ such that $N (u) \cap D \neq N(v) \cap D$ for every two distinct vertices $u$ and $v$ in $V(G) \backslash D$.

The regular independence number was first introduced by Albertson and Boutin in \cite{AB}. They proved lower bounds for planar graphs, maximal planar graphs, bounded-degree graphs and trees. Recently, Caro, Hansberg and Pepper \cite{CHP} generalised the regular independence number by introducing the regular $k$-independence number $\alpha_{k-reg}(G)$ of a graph $G$, and they generalised the results in \cite{AB} and found lower bounds for the regular $k$-independence numbers of trees, forests, planar graphs, $k$-trees and $k$-degenerate graphs. Guo, Zhao, Lai and Mao \cite{GZLM} obtained the exact values of the regular $k$-independence numbers of some special classes of graphs, and they established some lower bounds and upper bounds for line graphs and trees with a given diameter. They also obtained results of Nordhaus--Gaddum \cite{NG} type.

The regular domination number was first introduced and studied by Caro, Hansberg and Henning \cite{CHH}. They referred to the regular domination number as the fair domination number. Das and Desormeaux \cite{DD} considered the problem of minimizing the size of a regular dominating set that induces a connected subgraph. Further results on fair domination are obtained in \cite{MIC, CHHM}.

For standard definitions and notation in graph theory, we refer to \cite{West}. For a graph $G$ and a subset $A$ of $V(G)$, $E(A,V(G) \backslash A)$ denotes the set of edges of $G$ which have one vertex in $A$ and the other in $V(G) \backslash A$. Unless specified otherwise, we make use of the following notation: $n=|V(G)|$, $m=|E(G)|$, $e(A,V(G) \backslash A)=|E(A,V(G) \backslash A)| $, $d(v) = |N(v)|$ and $\delta = \delta(G)$, $\Delta = \Delta (G)$. The \emph{maximum cut of $G$}, denoted by $\beta = \beta(G)$, is $\max \{ e(A, V(G) \backslash A) \colon A \subseteq V(G)\}$. For a non-negative integer $k$, we denote $\{i \colon 1 \leq i \leq k, \, i \mbox{ is an integer}\}$ by $[k]$. Note that $[k] = \emptyset$ if $k = 0$.

This paper is organized as follows. In Section~\ref{irrindsection}, we prove several sharp upper bounds for $\alpha_{ir}(G)$. In Section~\ref{irrind1section}, we characterize the graphs $G$ with $\alpha_{ir}(G)=1$, we determine those that are planar, and we determine those that are outerplanar. In Section~\ref{irrdomsection}, we prove several sharp lower bounds for $\gamma_{ir}(G)$, we characterize the graphs $G$ with $\gamma_{ir}(G) \in \{n, n-1\}$, and we also provide some upper bounds for $\gamma_{ir}(G)$. In Section~\ref{irrinddomsection}, we provide sharp upper bounds relating $\alpha_{ir}(G)$ to $\gamma_{ir}(G)$ or $\gamma_{ir}(\bar{G})$. In Section~\ref{NDsection}, we provide sharp Nordhaus--Gaddum-type bounds for both  $\alpha_{ir}(G)$ and $\gamma_{ir}(G)$.

\section{Irregular independence} \label{irrindsection}

In this section, we provide various bounds for $\alpha_{ir}(G)$. We start with bounds in terms of basic graph parameters.

For any graph $G$, we denote by $span(G)$ the number of distinct values in the degree sequence of $G$. More formally, $span(G)=| \{ d(v) \colon v \in V(G) \}|$. Clearly, $span(G) \leq \Delta - \delta +1$.

\begin{theorem-section} \label{irregind1}
For any graph $G$,
%
\[1 \leq \alpha_{ir}(G) \leq \min \left\{ \Delta - \delta +1, \, \left \lfloor \frac{n-\delta+1}{2} \right \rfloor , \, \frac{1+ \sqrt{2n^2-2n-4m+1}}{2} \right\}.\]
%
Moreover, the bound is sharp.

\end{theorem-section}

\begin{proof}
We have $\alpha_{ir} (G) \geq 1$ as $\{ v \}$ is an irregular independent set for each $v \in V(G)$. Clearly, $\alpha _ {ir} (G) \leq span(G) \leq \Delta - \delta +1$. Let $A$ be a largest irregular independent set. Let $v_1, \dots, v_t$ be the distinct vertices of $A$ with $\delta \leq d (v_1) < \dots < d (v_t)$. Thus, $\delta +t-1 \leq d(v_t) \leq |V(G) \backslash A| =n-t$, from which we get $t \leq \left \lfloor \frac{n-\delta+1}{2} \right \rfloor$. Let $B = V(G) \backslash A$. We have 
\begin{align} m &= |E(G[B])| + \sum_{v \in A} d(v) \leq \frac{1}{2}(n-t)(n-t-1) + \sum _{i=1} ^{t} (n-2t + i) \nonumber \\
&= \frac{1}{2}(n-t)(n-t-1) + \frac{t}{2} (2n-3t+1), \nonumber
\end{align}
so $2t^2-2t + (n+2m-n^2) \leq 0$, and hence $\alpha _{ir} (G) \leq \frac{1}{2}\left(1+ \sqrt{2n^2-2n-4m+1}\right)$. This establishes the bound in the theorem.

The lower bound is attained if $G$ is regular. We now show that the upper bound is sharp. Let $r$ and $t$ be positive integers. 

If $G$ is the disjoint union of $K_r, K_{r+1}, \dots, K_{r+t-1}$, then $\alpha_{ir}(G) = \Delta - \delta + 1$. 

Let $k = r + t - 1$. Suppose that $G$ is constructed as follows: let $v_1, \dots, v_t, w_1, \dots, w_k$ be the distinct vertices of $G$, and, for each $i \in [t]$, form exactly $r +i-1$ distinct edges of the form $\{v_i, w_j\}$. Let $A = \{v_1, \dots, v_t\}$ and $B = \{w_1, \dots, w_k\}$. Since $A$ is an irregular independent set of $G$, $\alpha_{ir}(G) \geq t$. But $\alpha_{ir}(G) \leq \lfloor \frac{n-\delta+1}{2} \rfloor = \lfloor \frac{(\delta+2t-1)-\delta+1}{2} \rfloor =t$. Thus, $\alpha_{ir}(G) = \lfloor \frac{n-\delta+1}{2} \rfloor$. 

Let $r \geq t$. Suppose that $G$ is constructed as follows: let $v_1, \dots, v_t, w_1, \dots, w_r$ be the distinct vertices of $G$, form a complete graph on the vertices $w_1, \dots, w_r$, and, for each $i \in [t]$, form exactly $r-t+i$ distinct edges of the form $\{v_i, w_j\}$. Let $A = \{v_1, \dots, v_t\}$. Since $A$ is an irregular independent set of $G$, $t \leq \alpha_{ir}(G)$. We have $m = \frac{1}{2}r(r-1) + \sum_{i=1}^{t} (r-t+i) = \frac{1}{2}r(r-1) + \frac{1}{2}t(2r-t+1)$. Since $n = r+t$, $2m = (n-t)(n-t-1) + t(2n-3t+1) = n^2 - n - 2t^2  + 2t$. 
%
%
By the established bound, $\alpha_{ir}(G) \leq \frac{1}{2}\left(1+ \sqrt{2n^2-2n-4m+1}\right) \leq t$. Since $\alpha_{ir}(G) \geq t$, $\alpha_{ir}(G) = \frac{1}{2}\left(1+ \sqrt{2n^2-2n-4m+1}\right)$.
\end{proof}

We also have
\begin{equation} \alpha_{ir} (G) \leq \frac{-2 \delta +1 + \sqrt{{(2 \delta -1)}^2 +8m}}{2}. \label{irregind3}
\end{equation}
This is immediate from our next result, the proof of which also shows that (\ref{irregind3}) is sharp. 

\begin{theorem-section} \label{irregind2}
For any graph $G$,
%
\[\alpha _{ir} (G) \leq \frac{-2 \delta +1 + \sqrt{{(2 \delta -1)}^2 +8\beta}}{2}.\]
%
Moreover, the bound is sharp.
\end{theorem-section}

\begin{proof}
Let $t = \alpha_{ir}(G)$. Let $A$ be an irregular independent set of $G$ of size $t$, and let $v_1, \dots, v_t$ be the distinct vertices in $A$. We have $\beta \geq e(A, V(G) \backslash A) = \sum_{i = 1}^t d (v_i) \geq \sum _ {i=0} ^ {t-1} ( \delta + i )  = \frac{1}{2} t(2 \delta + t -1)$, so $0 \geq t^2 + (2\delta-1)t - 2 \beta$. Solving the quadratic inequality, we obtain $t \leq \frac{1}{2}\left(-2 \delta +1 + \sqrt{{(2 \delta -1)}^2 +8\beta} \right)$. 

We now prove that the bound is sharp. Let $r$ and $t$ be positive integers such that $t(t-1) \geq 2r(r-1)$. Let $k = r + t - 1$. Let mod$^*$ be the usual modulo operation with the exception that, for every two positive integers $a$ and $b$, $ba  \; {\rm mod}^* \; a$ is $a$ rather than $0$. Let $s_0 = 0$, and let $s_i = \sum_{j=0}^{i-1} (r+j)$ for each $i \in [t]$. Suppose that $G$ is constructed as follows: let $v_1, \dots, v_t, w_1, \dots, w_k$ be the distinct vertices of $G$, and, for each $i \in [t]$, let $v_i$ be adjacent to the vertices in $\{w_{j \; {\rm mod^*} \; k} \colon j \in [s_{i-1}+1, s_i]\}$. Thus, $v_1$ is adjacent to $w_1, \dots, w_r$, $v_2$ is adjacent to $w_{r+1 \; {\rm mod}^* \; k}, \dots, w_{2r+1 \; {\rm mod}^* \; k}$, $v_3$ is adjacent to $w_{2r+2 \; {\rm mod}^* \; k}, \dots, w_{3r+3 \; {\rm mod}^* \; k}$, and so on. By construction, $d(w_k) = \min\{d(w_j) \colon j \in [k]\}$. Let $A = \{v_1, \dots, v_t\}$ and $B = \{w_1, \dots, w_k\}$. Since $G$ is a bipartite graph with partite sets $A$ and $B$, we have $\beta = m = e(A,B) = \sum_{i = 1}^t d(v_t) = s_t = \frac{1}{2}t(2r + t-1)$. We also have $m = \sum_{j = 1}^k d(w_j) \geq d(w_k)k$, so $\frac{1}{2}t(2r + t-1) \leq d(w_k)k$ and hence $d(w_k) \geq \frac{t(2r+t-1)}{2k} = \frac{t(2r+t-1)}{2(r+t-1)}$. If we assume that $\frac{t(2r+t-1)}{2(r+t-1)} < r$, then we get a contradiction to the condition $t(t-1) \geq 2r(r-1)$. Thus, $d(w_k) \geq r$. Since $\min\{d(v_i) \colon i \in [t]\} = d(v_1) = r \leq d(w_k) = \min\{d(w_j) \colon j \in [k]\}$, $\delta = d(v_1) = r$. Now $A$ is an irregular independent set of $G$, so $\alpha_{ir}(G) \geq t$. By the bound in the theorem, 
\begin{align} \alpha_{ir}(G) &\leq \frac{-2 \delta +1 + \sqrt{{(2 \delta -1)}^2 +8\beta}}{2}  = \frac{-2r +1 + \sqrt{{(2r -1)}^2 + 4t(2r +t-1)}}{2} \nonumber \\
&= \frac{-2r +1 + \sqrt{(2r +2t-1)^2}}{2} = t. \nonumber
\end{align}
Since $\alpha_{ir}(G) \geq t$, $\alpha _{ir} (G) = \frac{1}{2}\left(-2 \delta +1 + \sqrt{{(2 \delta -1)}^2 +8\beta} \right)$.
\end{proof}


%
%


Our next result provides inequalities relating $\alpha_{ir}(G)$ to $\alpha_{reg}(G)$.

\begin{theorem-section} \label{ir-reg-ind}
For any graph $G$,
\begin{enumerate}[(i)]
\item $2 \leq \alpha _ {ir} (G) + \alpha _{reg} (G) \leq n+1$. 
\item $\alpha (G) \leq \alpha_{ir}(G)  \alpha_{reg}(G) \leq {(\alpha(G))}^2$.
\item if $n \geq 4$, then $1 \leq \alpha _ {ir} (G)  \alpha _{reg} (G) \leq \lfloor \frac{n}{2} \rfloor \lceil \frac{n}{2} \rceil$.
\end{enumerate}
Moreover, the following hold:
\begin{enumerate}[(a)]
\item The bounds are sharp.
\item The upper bound in (i) is attained if and only if $G$ is empty. Also, for any integer $k$ with $2 \leq k \leq n+1$, $\alpha_{ir}(G) + \alpha_{reg}(G) = k$ if $G = E_{k-2} \cup K_{n-k+2}$. 
\end{enumerate}
\end{theorem-section}

\begin{proof}
Let $A$ be an irregular independent set of $G$ of size $\alpha _ {ir} (G)$. Let $B$ be a regular independent set of $G$ of size $\alpha _ {reg} (G)$. Let $I$ be a largest independent set of $G$.

(i) Trivially, $\alpha_{ir}(G) \geq 1$, $\alpha_{reg}(G) \geq 1$, and hence the lower bound.  Clearly, $|A \cap B| \leq 1$. We have $n \geq |A \cup B| = |A| + |B| - |A \cap B| \geq \alpha _{ir}(G) + \alpha _ {reg} (G) -1$, so $\alpha _ {ir} (G) + \alpha _{reg} (G) \leq n+1$. 

(ii) Let $d_1,\dots,d_r$ be the distinct degrees of the vertices in $I$. For each $i \in [r]$, let $D_i$ be the set of vertices in $I$ of degree $d_i$. Let $s = \max\{|D_i| \colon i \in [r] \}$. We have $r \leq \alpha_{ir}(G)$, $s \leq \alpha_{reg}(G)$, and $\alpha(G) = |I|= |D_1| + \dots +|D_r| \leq rs \leq \alpha _{ir} (G) \alpha_{reg} (G)$. Trivially, $\alpha_{ir}(G) \leq \alpha (G)$, $\alpha_{reg} (G) \leq \alpha (G)$, and hence the upper bound.

(iii) As in (i), the lower bound is trivial. By (i), $|A| + |B| \leq n+1$. Suppose equality holds. Then $G = E_n$ by (b), which is proved below. Thus, $|A| |B| = n \leq \lceil \frac{n}{2} \rceil \lfloor \frac{n}{2} \rfloor$ if $n \geq 4$. Now suppose $|A| + |B| \leq n$. Then $|A||B| \leq |A|(n-|A|)$. By differentiating the function $f(r) = r(n-r)$, we see that $f$ increases as $r$ increases from $0$ to $\frac{n}{2}$. Thus, $|A||B| \leq \lfloor \frac{n}{2} \rfloor (n - \lfloor \frac{n}{2} \rfloor ) = \lfloor \frac{n}{2} \rfloor \lceil \frac{n}{2} \rceil$. Hence the upper bound.\medskip

(a) The lower bounds in (i)--(iii) and the upper bound in (ii) are attained if $G = K_n$. The upper bound in (i) is attained if $G = E_n$. 

We now show that the upper bound in (iii) is sharp. For each of Cases 1--4 below, we construct a graph that attains the bound. Let $v_1, \dots, v_n$ be its distinct vertices. If $n \bmod 4 = 0$, then let $X = \{v_1, \dots, v_{\frac{n}{2}}\}$, let $Y = \{v_{\frac{n}{2} + 1}, \dots, v_n\}$, and, for each $j \in [n/4]$, let $v_j$ be adjacent to exactly $j-1$ vertices in $Y$, and let $v_{\frac{n}{2}-j+1}$ be adjacent to the remaining vertices in $Y$. If $n \bmod 4 = 1$, then let $X = \{v_1, \dots, v_{\frac{n-1}{2}}\}$, let $Y = \{v_{\frac{n+1}{2} + 1}, \dots, v_n\}$, and, for each $j \in [(n-1)/4]$, let $v_j$ be adjacent to exactly $j$ vertices in $Y$, and let $v_{\frac{n-1}{2}-j+1}$ be adjacent to the remaining vertices in $Y$. If $n \bmod 4 = 2$, then let $X = \{v_1, \dots, v_{\frac{n}{2}}\}$, let $Y = \{v_{\frac{n}{2}+1}, \dots, v_n\}$, let $v_{\frac{n}{2}}$ be adjacent to each vertex in $Y$, and, for each $j \in [(n-2)/4]$, let $v_j$ be adjacent to exactly $j$ vertices in $Y$, and let $v_{\frac{n}{2}-j}$ be adjacent to the remaining vertices in $Y$. If $n \bmod 4 = 3$, then let $X = \{v_1, \dots, v_{\frac{n+1}{2}}\}$, let $Y = \{v_{\frac{n+3}{2}}, \dots, v_n\}$, and, for each $j \in [(n+1)/4]$, let $v_j$ be adjacent to exactly $j-1$ vertices in $Y$, and let $v_{\frac{n+1}{2}-j+1}$ be adjacent to the remaining vertices in $Y$. Suppose that the resulting graph is $G$. Then $X$ is an irregular independent set of $G$, $Y$ is a regular independent set of $G$, and $|X| |Y| = \lfloor \frac{n}{2} \rfloor \lceil \frac{n}{2} \rceil$. By the bound in (iii), $\alpha _ {ir} (G)  \alpha _{reg} (G) = \lfloor \frac{n}{2} \rfloor \lceil \frac{n}{2} \rceil$.\medskip

(b) As stated in (a), the upper bound in (i) is attained in $G = E_n$. We now prove the converse. Thus, suppose $\alpha _ {ir} (G) + \alpha _{reg} (G) = n+1$. Thus, $|A| + |B| = n+1$. Recall that $|A \cap B| \leq 1$. Thus, $n \leq |A| + |B| - |A \cap B| = |A \cup B| \leq n$, giving $|A \cup B| = n$ and $|A \cap B| = 1$. Thus, for some $v \in V(G)$, $A \cap B = \{v\}$ and $A = (V(G) \backslash B) \cup \{v\}$. If $d(v)=0$, then since $v \in B$, all the vertices of $B$ must have degree $0$. Since $A$ and $B$ are independent sets containing $v$, $v$ has no neighbours in $A \cup B$. Thus, $d(v) = 0$ as $A \cup B = V(G)$. Hence $d(w) = 0$ for each $w \in B$. Now consider any $x \in V(G) \backslash B$. We have $x \in A$. Since $A$ is independent, $N(x) \subseteq B$. Since the vertices in $B$ have no neighbours, $N(x) = \emptyset$. Thus, $G$ is empty, as required. 

It is easy to check that $\alpha_{ir}(G) + \alpha_{reg}(G) = k$ if $G = E_{k-2} \cup K_{n-k+2}$ with $2 \leq k \leq n+1$.
\end{proof} 

\begin{corollary-section}
For any graph $G$ on $n \geq 4$ vertices,
\begin{center}
$\alpha_{ir} (G) \times \alpha_{reg} (G) \leq \min \{(\alpha (G))^2 , \lceil \frac{n}{2} \rceil \lfloor \frac{n}{2} \rfloor \}$.
\end{center}
\end{corollary-section}

\section{Graphs with irregular independence number $1$} \label{irrind1section}

We now investigate the particularly interesting case $\alpha_{ir}(G) = 1$. Let $D(G)$ denote the set of degrees of vertices of $G$. For any $i \in D(G)$, let $N_i$ denote the set of vertices of $G$ of degree $i$. Let $n_i = |N_i|$. For any two disjoint subsets $X$ and $Y$ of $V(G)$, let ${<}X,Y{>}$ denote the subgraph of $G$ given by $(X \cup Y , \{ \{x,y \} \in E(G) \colon x \in X, y \in Y \})$. 

\begin{lemma-section} \label{lemmairind1}
If $\alpha_{ir}(G) = 1$, then 
\begin{enumerate}[(i)]
\item ${<}N_i, N_j{>}$ is a complete bipartite graph for any $i,j \in D(G)$ with $i \neq j$.
\item the subgraph of $G$ induced by $N_k$ is $(k+n_k-n)$-regular for any $k \in D(G)$.
\end{enumerate}
\end{lemma-section}

\begin{proof}
(i) Suppose $\{v,w\} \notin E(G)$ for some $v \in N_i$ and some $w \in N_j$ with $i \neq j$. Then $\{v, w\}$ is an irregular independent set of $G$ of size $2$. This contradicts $\alpha_{ir}(G) = 1$. 

(ii) Let $v \in N_k$. By (i), for any $j \in D(G) \backslash \{k\}$, $v$ is adjacent to each $w \in N_j$. Thus, $v$ is adjacent to each vertex in $V(G) \backslash N_k$. By definition of $N_k$, the degree of $v$ in the subgraph of $G$ induced by $N_k$ is $k - (n - n_k)$.
\end{proof}

\begin{theorem-section} \label{theoremirind1}
If $\alpha_{ir}(G)=1$, then
\begin{enumerate}[(i)]
\item $n_k \geq n-k$ for any $k \in D(G)$.
\item $span(G) \leq \frac{1}{2}(1+\sqrt{1+8 \delta})$. Moreover, the bound is sharp.
\end{enumerate}

\end{theorem-section}

\begin{proof} (i) By Lemma~\ref{lemmairind1}(ii), $k + n_k - n \geq 0$.

(ii) Let $t = span(G)$. Then $D(G) = \{d_1, \dots, d_t\}$  for some integers $d_1, \dots, d_t$ with $0 \leq d_1 < \dots < d_t$. Now 
\begin{align} n &= \sum_{i=1}^{t} n_{d_i} \geq \sum_{i=1}^{t} (n-d_i) = (n-d_1) + \sum_{i=2}^{t} (n-d_i) = (n- \delta) + \sum_{i=2}^{t} n - \sum_{i=2}^{t} d_i \nonumber \\
&\geq (n- \delta) +  (t-1)n - \sum_{i=2}^{t} (n-t+i-1)= tn- \delta - \frac{(t-1)}{2} (2n-t). \nonumber
\end{align}
Therefore, $0 \geq t^2 -t-2\delta$, and the bound follows. The bound is attained if, for example, $G$ is the complete $k$-partite graph $K_{1,\dots,k}$. Indeed, we then have $\alpha_{ir}(G) =1$, $\delta = n-k$, $n=1+\dots+k = \frac{k}{2}(k+1)$, and $k = span(G) \leq \frac{1+ \sqrt{1+8 \delta}}{2} = \frac{1+ \sqrt{1+8(n-k)}}{2} = \frac{1+ \sqrt{1+8(\frac{k}{2}(k+1)) -8k}}{2} = \frac{1+\sqrt{(2k-1)^2}}{2} = k$, so $span(G) = \frac{1+ \sqrt{1+8 \delta}}{2}$.
\end{proof}


Let $G$ and $H$ be two vertex-disjoint graphs. The \emph{join of $G$ and $H$}, denoted by $G+H$, is the graph with $V(G+H)=V(G) \cup V(H)$ and $E(G+H)=E(G) \cup E(H) \cup \{ \{ x,y\} \colon x \in V(G), y\in V(H) \}$. If $k \geq 2$, $r \geq 2$, $G = K_1$, and $H$ is a vertex-disjoint union of $r$ copies of $K_{k-1}$, then $G+H$ is called a \emph{$k$-windmill graph} and is denoted by $Wd(k,r)$.

\begin{theorem-section} \label{planarresult}
A graph $G$ is planar and $\alpha_{ir}(G)= 1$ if and only if $G$ is a regular planar graph or a copy of one of the graphs $K_{1,n-1}$, $K_{2,n-2}$, $K_2 + E_{n-2}$, $K_2 + \frac{n-2}{2} K_2$, $E_2 + \frac{n-2}{2} K_2$, $E_2 + C_{n-2}$, $Wd(3,\frac{n-1}{2})$, and $K_1+H$, where $H$ is a union of vertex-disjoint cycles.
\end{theorem-section}

We start the proof of the theorem above with the following two lemmas.

\begin{lemma-section}\label{outerlemma} If a planar graph $G$ has a vertex $v$ that is adjacent to all the other vertices of $G$, then $G-v$ is outerplanar.
\end{lemma-section}

\begin{proof} Indeed, by deleting $v$ (and all edges incident to it) from a plane drawing of $G$, we obtain a  plane drawing of $G-v$ that has all the vertices on the same face. This means that $G-v$ is outerplanar because, for any face $F$ of a plane drawing $\varphi$ of a planar graph, $\varphi$ can be transformed to another plane drawing of the same graph in such a way that $F$ becomes the unbounded face, for example, by using stereographic projection (see \cite[Remark~6.1.27]{West}).
\end{proof}

\begin{lemma-section}\label{planarjoinlemma} If $\varphi$ is a plane drawing of $E_2 + C_k$ ($k \geq 3$), then a vertex $v$ of $E_2$ is mapped by $\varphi$ into the interior $I$ of the drawing of $C_{k}$, and the other vertex $w$ of $E_2$ is mapped by $\varphi$ into the exterior $E$ of the drawing of $C_{k}$.
\end{lemma-section}

\begin{proof} Let $G = E_2 + C_k$. Let $F \in \{I,E\}$ such that $v$ is mapped by $\varphi$ into $F$. Since $v$ is adjacent to each vertex of $C_k$, each face of $F$ in the drawing of $G - w$ has exactly $3$ vertices on its boundary, one of which is $v$. Thus, if we assume that $w$ is mapped into $F$, then we obtain that $w$ lies in the interior of one of these faces, and hence that $w$ is adjacent to at most two vertices of $C_k$, a contradiction.
\end{proof}

\begin{proof}[Proof of Theorem~\ref{planarresult}]
It is easy to check that if $G$ is one of the explicit graphs in Theorem~\ref{planarresult}, then $G$ is planar and $\alpha_{ir}(G)=1$. We now prove the converse. 

Let $G$ be a planar graph with $\alpha_{ir}(G)=1$. Since $K_5$ and $K_{3,3}$ are non-planar, $G$ does not contain any copies of these. It is well known that having $G$ planar implies that $m \leq 3n-6$. 
%
%
Suppose that $G$ is not regular. Setting $t = span(G)$, we then have $t \geq 2$ (and $n \geq 3$). We have $D(G) = \{d_1, \dots, d_t\}$  for some integers $d_1, \dots, d_t$ with $0 \leq d_1 < \dots < d_t$. We will often use Lemma~\ref{lemmairind1}(i), which tells us that, for any $i,j \in D(G)$ with $i \neq j$, each vertex of $N_{d_i}$ is adjacent to each vertex of $N_{d_j}$. The first immediate deduction from this is that $d_1 \geq 1$ as $t \geq 2$.

Suppose $t \geq 3$. Let $\{a_1, \dots, a_t\} = \{d_1, \dots, d_t\}$ such that $n_{a_1} \leq \dots \leq n_{a_t}$. If we assume that $n_{a_1} = n_{a_2} = 1$, then Lemma~\ref{lemmairind1}(i) gives us $a_1 = a_2 = n-1$, a contradiction (as $a_1, \dots, a_t$ are distinct). Thus, $n_{a_i} \geq 2$ for each $i \in [2,t]$. If we assume that $\sum_{i=3}^t n_{a_i} \geq 3$, then, by Lemma~\ref{lemmairind1}(i), we obtain that ${<}N_{a_1} \cup N_{a_2}, \bigcup_{i=3}^t N_{a_i}{>}$ contains a copy of $K_{3,3}$, a contradiction. Thus, $t = 3$ and $n_{a_2} = n_{a_3} = 2$. Let $\{u_1, u_2\} = N_{a_2}$ and $\{v_1, v_2\} = N_{a_3}$. We cannot have $\{u_1, u_2\}, \{v_1, v_2\} \in E(G)$, because otherwise Lemma~\ref{lemmairind1}(i) gives us $a_2 = n_{a_1} + n_{a_3} + 1 = n_{a_1} + 3 = n_{a_1} + n_{a_2} + 1 = a_3$, a contradiction. Similarly, we cannot have $\{u_1, u_2\}, \{v_1, v_2\} \notin E(G)$. Thus, for some $i \in \{2,3\}$, $a_i = n_{a_1} + 2$ and $a_{5-i} = n_{a_1} + 3$. We cannot have $n_{a_1} = 1$, because otherwise $a_1 = n_{a_2} + n_{a_3} = 4 = a_{5-i}$. Thus, $n_{a_1} = 2$. Let $\{w_1,w_2\} = N_{a_1}$. We cannot have $\{w_1,w_2\} \in E(G)$, because otherwise $a_1 = 5 = a_{5-i}$. Thus, we have $\{w_1,w_2\} \notin E(G)$, which gives us $a_1 = 4 = a_i$, a contradiction.  

Therefore, $t=2$. If we assume that $n_{d_1} \geq 3$ and $n_{d_2} \geq 3$, then, by Lemma \ref{lemmairind1}(i), we obtain that $G$ contains a copy of $K_{3,3}$, a contradiction. Thus, $n_{d_i} \leq 2$ for some $i \in \{1,2\}$. Let $j = 3-i$. By Lemma~\ref{lemmairind1}(i), $G = G[N_{d_i}] + G[N_{d_j}]$. By Lemma~\ref{lemmairind1}(ii), $G[N_{d_j}]$ is $k$-regular, where $k = d_j + n_{d_j} - n$.

Suppose $n_{d_i} = 1$. Let $\{v\} = N_{d_i}$. 
Thus, $G = (\{v\},\emptyset) + G[N_{d_j}]$. By Lemma~\ref{outerlemma}, $G[N_{d_j}]$ is outerplanar. Since the minimum degree of an outerplanar graph is at most $2$ (see \cite[Proposition~6.1.20]{West}), $k \leq 2$. 
If $k=0$, then $G$ is a copy of $ K_{1,n-1}$. If $k=1$, then $G[N_{d_j}] $ is a copy of $\frac{n-1}{2}K_2$, so $G$ is a copy of $Wd(3,\frac{n-1}{2})$. If $k=2$, then $G[N_{d_j}]$ is a cycle or a union of vertex-disjoint cycles.

Now suppose $n_{d_i} = 2$. Let $\{v,w\} = N_{d_i}$ and let $\{u_1 ,\dots , u_{n-2}\} = N_{d_j}$. By the handshaking lemma, $|E(G[N_{d_j}])| = \frac{k(n-2)}{2}$. By Lemma~\ref{lemmairind1}(i), $|E({<}N_{d_i}, N_{d_j}{>})| = 2(n-2)$. Now $m = |E(G[N_{d_i}])| + |E(G[N_{d_j}])| + |E({<}N_{d_i}, N_{d_j}{>})| \geq \frac{k(n-2)}{2} + 2(n-2) $. Since $m \leq 3n-6$, we obtain $k \leq 2$.

If $k=0$ and $\{v, w\} \in E(G)$, then $G$ is a copy of $K_2 + E_{n-2}$. If $k=0$ and $\{v, w\} \notin E(G)$, then $G$ is a copy of $E_2 + E_{n-2} = K_{2,n-2}$. If $k=1$ and $\{v, w\} \in E(G)$, then $G$ is a copy of $K_2 + \frac{n-2}{2}K_2$. If $k=1$ and $\{v, w\} \notin E(G)$, then $G$ is a copy of $E_2 + \frac{n-2}{2}K_2$. 

Finally, suppose $k=2$. We cannot have $v$ adjacent to $w$, because otherwise $m = 1 + \frac{2(n-2)}{2} + 2(n-2) > 3n-6$. Since $k=2$, $G[N_{d_j}]$ is a union of vertex-disjoint cycles $G_1, \dots, G_r$. Suppose $r \geq 2$. Let $\theta$ be a plane drawing of $G$. Let $\varphi$ be the drawing obtained by restricting $\theta$ to the subgraph $G' = (\{v,w\}, \emptyset) + G_1$ of $G$. By Lemma~\ref{planarjoinlemma}, no face of $\varphi$ has both $v$ and $w$ on its boundary. Since $G'$ and $G_2$ are vertex-disjoint, the drawing of $G_2$ in $\theta$ lies in the interior of one of the faces of $\varphi$. Thus, no vertex of $G_2$ is adjacent to both $v$ and $w$. This contradicts $G = G[N_{d_i}] + G[N_{d_j}]$. Therefore, $r = 1$. Thus, $G$ is $G[N_{d_i}] + G_1$, which is a copy of $E_2 + C_{n-2}$. 
\end{proof}

\begin{corollary-section}\label{outerplanarresult}
A graph $G$ is outerplanar and $\alpha_{ir}(G) = 1$ if and only if $G$ is a union of vertex-disjoint cycles or a copy of one of the graphs $E_n$, $\frac{n}{2}K_2$, $K_{1,n-1}$, $K_{2,2}$, $K_2 + E_2$, and $Wd(3,\frac{n-1}{2})$.
\end{corollary-section}

\begin{proof}
It is trivial that if $G$ is one of the explicit graphs in the statement of Corollary~\ref{outerplanarresult}, then $G$ is outerplanar and $\alpha_{ir}(G) = 1$.

We now prove the converse. Let $G$ be an outerplanar graph with $\alpha_{ir}(G)= 1$. This means that $\delta \leq 2$, as mentioned in the proof of Theorem~\ref{planarresult}. If $G$ is $k$-regular, then $k \leq 2$ and hence $G$ is a copy of $E_n$ (if $k = 0$) or a copy of $\frac{n}{2}K_2$ (if $k = 1$) or a union of vertex-disjoint cycles (if $k = 2$). Suppose that $G$ is not regular. Since $\delta \leq 2$, it follows by Theorem~\ref{planarresult} that $G$ is a copy of one of $K_{1,n-1}$, $K_{2,n-2}$, $K_2 + E_{n-2}$, $E_2 + \frac{n-2}{2} K_2$, and $Wd(3,\frac{n-1}{2})$. Now $K_{2,3}$ is not outerplanar. Thus, $K_{2,n-2}$ is outerplanar only if $n \leq 4$. Also, for $n \geq 5$, $K_2 + E_{n-2}$ is not outerplanar as it contains $K_{2,3}$. Similarly, $E_2 + \frac{n-2}{2} K_2$ is planar only if $\frac{n-2}{2} \leq 1$. Hence the result.
\end{proof}

\section{Irregular domination} \label{irrdomsection}

In this section, we provide bounds for the irregular domination number, $\gamma_{ir}(G)$, and investigate cases of particular importance, primarily cases where a bound is attained. 

We will start with lower bounds for $\gamma_{ir}(G)$.

\begin{theorem-section} \label{irregdom22}
For any graph $G$,
\[\gamma _ {ir} (G) \geq \max \left\{ \left \lceil \frac{n}{2} \right \rceil , \, n- \Delta \right\}.\]
Moreover, the bound is sharp.

\end{theorem-section}

\begin{proof}
Let $t = \gamma_{ir} (G)$. Let $D$ be an irregular dominating set of $G$ of size $t$. Let $v_1, \dots, v_{n-t}$ be the vertices in $V(G) \backslash D$. For each $i \in [n-t]$, let $w_i = |N(v_i) \cap D|$; since $D$ is a dominating set, $w_i \geq 1$. We may assume that $w_1 < \dots < w_{n-t}$. We have $t=|D| \geq w_{n-t} \geq n-t$, and hence $t \geq \left \lceil \frac{n}{2} \right \rceil$. Since $n-t \leq w_{n-t} \leq \Delta$, $t \geq n - \Delta$. 

We now show that the bound is sharp. Let $k = \left \lceil \frac{n}{2} \right \rceil$ and $n' = n-k$. Suppose that $G$ is constructed as follows: let $u_1, \dots, u_k, v_1, \dots, v_{n'}$ be the distinct vertices of $G$, and, for each $i \in [n']$, let $v_i$ be adjacent to exactly $i$ of the vertices $u_1, \dots, u_k$. Since $\max\{d(u_i) \colon i \in [k]\} \leq n' = d(v_{n'}) = \max\{d(v_i) \colon i \in [n']\}$, $\Delta = n'$. Clearly, $\{u_1, \dots, u_k\}$ is an irregular dominating set of $G$ of size $\left \lceil \frac{n}{2} \right \rceil = n - n' = n - \Delta$. \end{proof}

\begin{theorem-section} \label{irregdom2}
For any graph $G$, 
\[\gamma_{ir}(G) \geq n + \frac{1- \sqrt{1+8 \beta}}{2}.\]
Moreover, the bound is sharp.
\end{theorem-section}

\begin{proof}
Let $t$, $D$, $v_1, \dots, v_{n-t}$, $w_1, \dots, w_{n-t}$ be as in the proof of Theorem~\ref{irregdom22}. 
We have $\beta \geq e(D,V(G) \backslash D) = \sum_{i=1}^{n-t} w_i \geq \sum_{i=1}^{n-t} i = \frac{1}{2}(n-t)(n-t+1)$, so $0 \geq t^2 - (2n+1)t + (n^2+n-2 \beta)$ and hence $t \geq n + \frac{1}{2}(1- \sqrt{1+8 \beta})$. 
 
We now show that the bound is sharp. Let $n/2 \leq k \leq n-1$ and $n' = n-k$. Suppose that $G$ is constructed as follows: let $u_1, \dots, u_k, v_1, \dots, v_{n'}$ be the distinct vertices of $G$, and, for each $i \in [n']$, let $v_i$ be adjacent to exactly $i$ of the vertices $u_1, \dots, u_k$. Let $D = \{u_1, \dots, u_k\}$. Since $D$ is an irregular dominating set of $G$, $\gamma_{ir}(G) \leq k$. Since $m  = e(D,V(G) \backslash D)$, $\beta = e(D,V(G) \backslash D) = \frac{1}{2}(n')(n'+1)$. By the established bound, 
\[\gamma_{ir}(G) \geq n + \frac{1- \sqrt{1+8\beta}}{2} = n + \frac{1- \sqrt{1+4(n')(n'+1)}}{2} = n + \frac{1- \sqrt{(2n-2k+1)^2}}{2} = k.\]
Since $\gamma_{ir}(G) \leq k$, $\gamma_{ir}(G) = n + \frac{1- \sqrt{1+8\beta}}{2}$. \end{proof}

\begin{corollary-section} If $G$ is a graph with average degree $d$, then  
\[\gamma_{ir}(G) \geq n - \sqrt{dn}.\] 
Moreover, equality holds if and only if $G$ is empty.
\end{corollary-section}

\begin{proof} Since $\beta \leq m$, $\gamma_{ir}(G) \geq n + \frac{1}{2}(1- \sqrt{1+8m})$ by Theorem~\ref{irregdom2}. Now $dn= \sum _ {v \in V(G)} d(v) = 2m$ (by the handshaking lemma), so $4dn = 8m$. Thus, $\gamma_{ir}(G)	\geq n + \frac{1}{2}(1- \sqrt{1+4dn}) \geq n + \frac{1}{2}(-\sqrt{4dn}) = n - \sqrt{dn}$. Note that equality holds throughout only if $d = 0$, in which case $G$ is empty.

If $G$ is empty, then $d = 0$ and $\gamma_{ir}(G) = n = n - \sqrt{dn}$.
\end{proof}

Next, we give a full characterization of the cases $\gamma _{ir} (G) = n$ and $\gamma _{ir} (G) = n-1$. For two graphs $G$ and $H$, we write $G \simeq H$ if $G$ is a copy of $H$. 

\begin{theorem-section} \label{irregdomemptyandregular}
For any graph $G$,
\begin{enumerate}[(i)] 
\item $\gamma_{ir}(G) = n$ if and only if $G \simeq E_n$.
\item $\gamma_{ir}(G) = n-1$ if and only if, for some $t \geq 0$ and some $r \geq 1$, $G \simeq tK_1 \cup K_{1,r}$ or $G \simeq tK_1 \cup H$ for some $r$-regular graph $H$.
\end{enumerate}
\end{theorem-section}

\begin{proof}
(i) If $G$ has an edge $\{ v,w \}$, then $V(G) \backslash \{ v \}$ is an irregular dominating set of $G$, so $\gamma _ {ir} (G) \leq n-1$. Therefore, $\gamma_{ir}(G) = n$ only if $G \simeq E_n$. If $G \simeq E_n$, then $V(G)$ is the only dominating set of $G$, so $\gamma_{ir}(G) = n$.

(ii) It is easy to see that $\gamma_{ir}(G) = n-1$ if $G \simeq tK_1 \cup K_{1,r}$ or $G \simeq tK_1 \cup H$ for some $r$-regular graph $H$. We now prove the converse. Thus, suppose $\gamma_{ir}(G) = n-1$. By (i), $E(G) \neq \emptyset$. 


Suppose $G$ has two vertices $u$ and $v$ such that $2 \leq d(u) < d(v)$. Then $V(G) \backslash \{u,v\}$ is an irregular dominating set of $G$ (independently of whether $u$ and $v$ are adjacent or not). Thus, we have $\gamma_{ir}(G) \leq n-2$, a contradiction. Therefore, 
\begin{equation} d(u) \leq 1 \mbox{ for any } u, v \in V(G) \mbox{ with } d(u) < d(v). \label{irregdomemptyandregular.1}
\end{equation}

Suppose $span(G) \geq 4$. Then there exist $v_1, v_2, v_3, v_4 \in V(G)$ such that $d(v_1) < d(v_2) < d(v_3) < d(v_4)$. Thus, we have $2 \leq d(v_3) < d(v_4)$, which contradicts (\ref{irregdomemptyandregular.1}). Therefore, $span(G) \leq 3$.

If $span(G) = 1$, then $G$ is an $r$-regular graph for some $r \geq 1$ ($r \neq 0$ as $E(G) \neq \emptyset$), and we are done. 

Suppose $span (G)=2$. Then $\{d(v) \colon v \in V(G)\} = \{p,r\}$ with $0 \leq p < r$. 
By (\ref{irregdomemptyandregular.1}), $p \leq 1$. If $p = 0$, then $G \simeq tK_1 \cup H$ for some $t \geq 1$ and some $r$-regular graph $H$. Suppose $p = 1$. Then $r \geq 2$. If we assume that there exists a pair of non-adjacent vertices $u$ and $v$ of degrees $1$ and $r$, respectively, then we obtain that $V(G) \backslash \{u, v\}$ is an irregular dominating set of $G$ of size $n-2$, which contradicts $\gamma_{ir}(G) = n-1$. Thus, each vertex $x$ of degree $1$ is adjacent to each vertex of degree $r$. Since $x$ has only one neighbour, there is only one vertex of degree $r$. Consequently, $G = K_{1,r}$.

Finally, suppose $span(G)=3$. Then there exist $v_1, v_2, v_3 \in V(G)$ such that $d(v_1) < d(v_2) < d(v_3)$. If we assume that $G$ has no vertex of degree $0$ or no vertex of degree $1$, then we obtain $2 \leq d(v_2) < d(v_3)$, which contradicts (\ref{irregdomemptyandregular.1}). Thus, since $span(G)=3$, $\{d(v) \colon v \in V(G) \} = \{0,1,r\}$ for some $r \geq 2$. Let $G'$ be the graph obtained by removing from $G$ the set $I$ of vertices of $G$ of degree $0$. Then $\{d(v) \colon v \in V(G')\} = \{1,r\}$. As in the case $span(G) = 2$ above, this yields $G' \simeq K_{1,r}$, so $G = tK_1 \cup K_{1,r}$, where $t = |I|$. \end{proof}

The \emph{Ramsey number} $R(p,q)$ is the smallest number $n$ such that every graph on $n$ vertices contains a clique of order $p$ or an independent set of order $q$.

\begin{theorem-section}
For any graph $G$,
\begin{enumerate}[(i)]
\item if $span(G) \geq R(k,k)$ and $\delta \geq k $, then $\gamma _ {ir} (G) \leq n-k$.
\item if $span(G) \geq 5$ and $\delta \geq 3$, then $\gamma_{ir} (G) \leq  n-3$.
\end{enumerate}

\end{theorem-section} 

\begin{proof}
(i) Suppose $span(G) \geq R(k,k)$ and $\delta \geq k$. Let $B$ be a set of $R(k,k)$ vertices of $G$ of distinct degrees. Then $G[B]$ has an independent set of size $k$ or a clique of size $k$. If $G[B]$ has an independent set $I$ of size $k$, then $V(G) \backslash I$ is an irregular dominating set of $G$ of size $n-k$. If $G[B]$ has a clique $K$ of size $k$, then, since $\delta \geq k$, $V(G) \backslash K$ is an irregular dominating set of $G$ of size $n-k$.

(ii) Suppose $span (G) \geq 5$ and $\delta \geq 3$. Let $B$ be a set of $5$ vertices of $G$ of distinct degrees. It is easy to see that if a $5$-vertex graph does not have an independent set of size $3$, then it is a copy of $C_5$ or has a clique of size $3$. If $G[B]$ is a copy of $C_5$, then each vertex in $B$ has a distinct number of neighbours in $V(G) \backslash B$, and hence, since $\delta \geq 3$, $V(G) \backslash B$ is an irregular dominating set of $G$ of size $n-5$. As in the proof of (i), $\gamma_{ir}(G) \leq n-3$ if $G[B]$ has an independent set of size $3$ or a clique of size $3$. 
\end{proof}

\section{Relations between irregular independence and irregular domination} \label{irrinddomsection}

We now establish a set of inequalities relating the irregular independence number to the irregular domination number. These are gathered in the theorem below. 
In the proof, we need to use the following more precise notation. For a vertex $v$ of a graph $G$, we will denote the set of neighbours of $v$ in $G$ by $N_G(v)$, and the degree of $v$ in $G$ by $d_G(v)$. Formally, $N_G(v) = \{w \in V(G) \colon vw \in E(G)\}$ and $d_G(v) = |N_G(v)|$. The \emph{complement of $G$} (that is, $(V(G), {V(G) \choose 2} \backslash E(G))$) is denoted by $\bar{G}$.

\begin{theorem-section}
For any graph $G$,
\begin{enumerate}[(i)]
\item $\alpha_{ir}(G) + \gamma_{ir}(G) \leq n+1$ if $\delta = 0$, and $\alpha_{ir}(G) + \gamma_{ir}(G) \leq n$ if $\delta \geq 1$.
\item $\alpha_{ir}(G) \gamma_{ir}(G) \leq \lfloor \frac{n+1}{2} \rfloor \lceil \frac{n+1}{2} \rceil$ if $\delta = 0$, and $\alpha_{ir}(G)  \gamma_{ir}(G) \leq \lfloor \frac{n}{2} \rfloor \lceil \frac{n}{2} \rceil$ if $\delta \geq 1$. 
\item $\alpha _ {ir} (G) + \gamma _{ir} (\bar{G}) \leq n+1$. 
\item $\alpha_{ir}(G)  \gamma_{ir}(\bar{G}) \leq \lfloor \frac{n+1}{2} \rfloor \lceil \frac{n+1}{2} \rceil$. 
\end{enumerate} 
Moreover, the bounds are sharp.

\end{theorem-section}

\begin{proof} Let $A$ be an irregular independent set of $G$ of size $\alpha_{ir}(G)$, and let $D = V(G) \backslash A$. 

Suppose $\delta \geq 1$. Then $D$ is an irregular dominating set of $G$, so $\alpha_{ir}(G) + \gamma_{ir}(G) \leq |A| + |D| \leq n$ and $\alpha_{ir}(G) \gamma_{ir}(G) \leq |A||D| = |A|(n-|A|) \leq \lfloor \frac{n}{2} \rfloor (n - \lfloor \frac{n}{2} \rfloor ) = \lfloor \frac{n}{2} \rfloor \lceil \frac{n}{2} \rceil$ (as in the proof of Theorem~\ref{ir-reg-ind}(iii)). Now suppose $\delta =0$. Let $V_0$ be the set of vertices of $G$ of degree $0$, and let $V_1$ be the set of vertices of $G$ of degree at least $1$. As in the case $\delta \geq 1$, $\alpha_{ir}(G[V_1]) + \gamma_{ir}(G[V_1]) \leq |V_1|$. We have $\alpha_{ir}(G) + \gamma_{ir}(G) = (\alpha_{ir}(G[V_1]) + 1) + (\gamma_{ir}(G[V_1]) + |V_0|) \leq  |V_0| + |V_1| + 1 = n+1$. Clearly, $A$ has exactly one element $x$ of $V_0$, and $D \cup \{x\}$ is an irregular dominating set of $G$. Thus, $\alpha_{ir}(G) \gamma_{ir}(G) \leq |A|(|D|+1) \leq |A|(n+1-|A|) \leq \lfloor \frac{n+1}{2} \rfloor (n+1 - \lfloor \frac{n+1}{2} \rfloor ) = \lfloor \frac{n+1}{2} \rfloor \lceil \frac{n+1}{2} \rceil$. Hence (i) and (ii).

Let $v_1, \dots, v_t$ be the distinct vertices in $A$, where $d_G(v_1) < \dots < d_G(v_t)$. We have $d_G(v_t) \leq |V(G) \backslash A| = n-t$. For each $i \in [t]$, let $a_i = |N_{\bar{G}}(v_i) \cap D|$. For each $i \in [t]$, $a_i = n - t - d_G(v_i) \geq n - t - d_G(v_t)$. Thus, if $d_G(v_t) \leq n-t-1$, then $D$ is an irregular dominating set of $\bar{G}$, and hence $\alpha_{ir}(G) + \gamma_{ir}(\bar{G}) \leq |A| + |D| = t + (n-t) = n$. Suppose $d_G(v_t) = n-t$. We have $a_i \geq 1$ for each $i \in [t-1]$. Let $A' = A \backslash \{v_t\}$. Let $D' = D \cup \{v_t\}$. For each $i \in [t-1]$, let $b_i = |N_{\bar{G}}(v_i) \cap D'|$. For each $i \in [t-1]$, we have $N_{\bar{G}}(v_i) \cap D' = (N_{\bar{G}}(v_i) \cap D) \cup \{v_t\}$, so $b_i = a_i + 1 = n - t - d_G(v_i) + 1$. Thus, $D'$ is an irregular dominating set of $\bar{G}$. Consequently, $\alpha_{ir}(G) + \gamma_{ir}(\bar{G}) \leq |A| + |D'| = t + (n-t+1) = n+1$ and $\alpha_{ir}(G)\gamma_{ir}(\bar{G}) \leq |A||D'| = t(n+1-t) \leq \lfloor \frac{n+1}{2} \rfloor (n+1 - \lfloor \frac{n+1}{2} \rfloor ) = \lfloor \frac{n+1}{2} \rfloor \lceil \frac{n+1}{2} \rceil$. Hence (iii) and (iv).\medskip

We now show that the bounds are sharp. We use constructions similar to that in the proof of Theorem~\ref{irregdom22}.

Let $k = \left \lceil \frac{n}{2} \right \rceil$ and $n' = n-k$. Suppose that $G$ is constructed as follows: let $u_1, \dots, u_k, v_1, \dots, v_{n'}$ be the distinct vertices of $G$, and, for each $i \in [n']$, let $v_i$ be adjacent to exactly $k-i+1$ of the vertices $u_1, \dots, u_k$. Clearly, $\delta \geq 1$. Also, $\{v_1, \dots, v_{n'}\}$ is an irregular independent set, and, by Theorem~\ref{irregind1}, it is of maximum size. Moreover, $\{u_1, \dots, u_k\}$ is an irregular dominating set of $G$, and, by Theorem~\ref{irregdom22}, it is of minimum size. Thus, $\alpha_{ir}(G) + \gamma_{ir}(G) = n' + k = n$ and $\alpha_{ir}(G) \gamma_{ir}(G) = n'k = \lfloor \frac{n}{2} \rfloor \lceil \frac{n}{2} \rceil$. Now suppose that we instead have that $k = \left \lceil \frac{n-1}{2} \right \rceil$, $n' = n-k$, and, for each $i \in [n']$, $v_i$ is adjacent to exactly $i-1$ of $u_1, \dots, u_k$. Since $d(v_1) = 0$, $\delta = 0$. Similarly to the above, $\{u_1, \dots, u_k, v_1\}$ is an irregular dominating set of $G$ of minimum size as $\{u_1, \dots, u_k\}$ is an irregular dominating set of $G-v_1$ of minimum size. Also, $\{v_1, \dots, v_{n'}\}$ is an irregular independent set of maximum size. Thus, $\alpha_{ir}(G) + \gamma_{ir}(G) = n' + k + 1 = n + 1$ and $\alpha_{ir}(G) \gamma_{ir}(G) = n'(k+1) = \lfloor \frac{n+1}{2} \rfloor \lceil \frac{n+1}{2} \rceil$. We have established that (i) and (ii) are sharp.

Let $k = \left \lceil \frac{n-1}{2} \right \rceil$ and $n' = n-k$. Suppose that $G$ is constructed as follows: let $u_1, \dots, u_k, v_1, \dots, v_{n'}$ be the distinct vertices of $G$, and, for each $i \in [n']$, let $v_i$ be adjacent to exactly $k-i+1$ of the vertices $u_1, \dots, u_k$. Thus, $\{v_1, \dots, v_{n'}\}$ is an irregular independent set, and, by Theorem~\ref{irregind1}, it is of maximum size (note that $\delta$ is $d(v_{n'})$, which is $0$ if $n$ is odd, and $1$ if $n$ is even). Also, we clearly have that $\{u_1, \dots, u_k, v_1\}$ is an irregular dominating set of $\bar{G}$, and it is of minimum size because $d_{\bar{G}}(v_1) = 0$ and, by Theorem~\ref{irregdom22}, $\{u_1, \dots, u_k\}$ is an irregular dominating set of $\bar{G}-v_1$ of minimum size. Thus, $\alpha_{ir}(G) + \gamma_{ir}(\bar{G}) = n' + k + 1 = n + 1$ and $\alpha_{ir}(G) \gamma_{ir}(\bar{G}) = n'(k+1) = \lfloor \frac{n+1}{2} \rfloor \lceil \frac{n+1}{2} \rceil$.
\end{proof}

\section{Nordhaus-Gaddum type results} \label{NDsection}

In this section, we provide results of Nordhaus--Gaddum type \cite{NG} for both the irregular independence number and the irregular domination number. We shall use the notation introduced in the preceding section. 
Also, where necessary, we will denote the minimum degree of $G$ and the maximum degree of $G$ by $\delta(G)$ and $\Delta(G)$, respectively. 

\begin{theorem-section}
If $G$ is a graph on $n \geq 2$ vertices, then

\begin{enumerate}[(i)] 
\item $2 \leq \alpha_{ir}(G) + \alpha_{ir}(\bar{G}) \leq n$, 
\item $1 \leq \alpha_{ir}(G) \alpha_{ir}(\bar{G}) \leq \lfloor \frac{n}{2} \rfloor \lfloor \frac{n+1}{2} \rfloor$. 
\end{enumerate}
Moreover, the bounds are sharp.
\end{theorem-section}

\begin{proof}
By Theorem~\ref{irregind1}, $1 \leq \alpha_{ir} (G) \leq \lfloor \frac{n- \delta (G) +1}{2} \rfloor$ and $1 \leq \alpha_{ir}(\bar{G}) \leq \lfloor \frac{n- \delta(\bar{G}) +1}{2} \rfloor$. The lower bounds follow immediately, and they are attained if $G$ is regular. If $\delta(G) \geq 1$, then $\alpha_{ir}(G) + \alpha_{ir}(\bar{G}) \leq \lfloor \frac{n}{2} \rfloor + \lfloor \frac{n + 1}{2} \rfloor \leq n$ and $\alpha_{ir}(G) \alpha_{ir}(\bar{G}) \leq \lfloor \frac{n}{2} \rfloor \lfloor \frac{n+1}{2} \rfloor$. Suppose $\delta(G) = 0$. Then $G$ has a vertex $v$ that has no neighbours. Thus, $v \in N_{\bar{G}}(u)$ for each $u \in V(\bar{G}) \backslash \{v\}$, and hence $\delta(\bar{G}) \geq 1$. This gives us $\alpha_{ir}(G) + \alpha_{ir}(\bar{G}) \leq \lfloor \frac{n+1}{2} \rfloor + \lfloor \frac{n}{2} \rfloor \leq n$ and $\alpha_{ir}(G) \alpha_{ir}(\bar{G}) \leq \lfloor \frac{n+1}{2} \rfloor \lfloor \frac{n}{2} \rfloor$.

We now show that the upper bounds are sharp. Let $k = \lceil \frac{n}{2} \rceil$ and $l = \lfloor \frac{n}{2} \rfloor$. Suppose that $G$ is constructed as follows: let $u_1, \dots, u_k, v_1, \dots, v_l$ be the distinct vertices of $G$, let every two distinct vertices in $\{v_1, \dots, v_l\}$ be adjacent, and, for each $i \in [k]$, let $u_i$ be adjacent to the vertices in $\{v_j \colon j \in [i-1]\}$. Clearly, $\{u_1, \dots, u_k\}$ is an irregular independent set of $G$, and $\{v_1, \dots, v_l\}$ is an irregular independent set of $\bar{G}$. Therefore, $\alpha_{ir}(G) + \alpha_{ir}(\bar{G}) \geq k + l = n$ and $\alpha_{ir}(G) \alpha_{ir}(\bar{G}) \geq kl$. By (i) and (ii), we actually have $\alpha_{ir}(G) + \alpha_{ir}(\bar{G}) = n$ and $\alpha_{ir}(G) \alpha_{ir}(\bar{G}) = kl$. Finally, note that $k = \lfloor \frac{n+1}{2} \rfloor$.
\end{proof}

\begin{theorem-section} \label{nordhausgaddum1}
If $G$ is a graph on $n \geq 2$ vertices, then

\begin{enumerate}[(i)] 
\item $2 \lceil \frac{n}{2} \rceil \leq \gamma _ {ir} (G) + \gamma _{ir} (\bar{G}) \leq 2n-1$, 
\item $( \lceil \frac{n}{2} \rceil )^2 \leq \gamma_{ir}(G) \gamma_{ir}(\bar{G}) \leq n(n-1)$. 
\end{enumerate}
Moreover, the following hold:
\begin{enumerate}[(a)]
\item The bounds are attainable for any $n \geq 3$.
\item For each of (i) and (ii), the upper bound is attained if and only if $G$ is empty or complete. 
\end{enumerate}
\end{theorem-section}

\begin{proof}
By Theorem \ref{irregdom22}, $\gamma_{ir}(G) \geq \left \lceil \frac{n}{2} \right \rceil$ and $\gamma_{ir}(G) \geq \left \lceil \frac{n}{2} \right \rceil$. The lower bounds in (i) and (ii) follow immediately. If $G$ is empty, then $\bar{G}$ is complete, so $\gamma_{ir}(G) + \gamma_{ir} (\bar{G}) = n + n-1 = 2n-1$ and $\gamma_{ir}(G) \gamma_{ir}(\bar{G}) = n(n-1)$. If $G$ is complete, then $\bar{G}$ is empty, so $\gamma_{ir}(G) + \gamma_{ir} (\bar{G}) = 2n-1$ and $\gamma_{ir}(G) \gamma_{ir}(\bar{G}) = (n-1)n$. If $G$ is neither empty nor complete, then $\bar{G}$ is non-empty, and hence, by Theorem~\ref{irregdomemptyandregular}, $\gamma_{ir} (G) + \gamma_{ir} (\bar{G}) \leq 2(n-1) < 2n-1$ and $\gamma_{ir}(G) \gamma_{ir}(\bar{G}) \leq (n-1)^2 < n(n-1)$. 

It remains to show that the lower bounds in (i) and (ii) are attainable for any $n \geq 3$. 

Suppose first that $n$ is odd. Let $k=\frac{n-1}{2}$. Suppose that $G$ is constructed as follows: let $u_1, \dots, u_k, v_1, \dots, v_{k+1}$ be the distinct vertices of $G$, and, for each $i \in [k]$, let $u_i$ be adjacent to $v_1, \dots, v_i$. Clearly, $\{v_1, \dots, v_{k+1}\}$ is an irregular dominating set of $G$ and of $\bar{G}$. Thus, $\gamma_{ir}(G) + \gamma_{ir}(\bar{G}) \geq 2(k+1) = 2 \lceil \frac{n}{2} \rceil$ and $\gamma_{ir}(G) \gamma_{ir}(\bar{G}) \geq (k+1)^2 = \lceil \frac{n}{2} \rceil^2$. By (i) and (ii), we actually have $\gamma_{ir}(G) + \gamma_{ir}(\bar{G}) = 2 \lceil \frac{n}{2} \rceil$ and $\gamma_{ir}(G) \gamma_{ir}(\bar{G}) = \lceil \frac{n}{2} \rceil^2$.

Now suppose that $n$ is even and $n \geq 8$. Let $k = \frac{n}{2}$. Suppose that $V(G) = \{u_1, \dots, u_k, v_1, \dots, v_k\}$ and that, for each $i \in [k] \backslash \{2\}$, $u_i$ is adjacent to $v_1, \dots, v_i$, $u_2$ is adjacent to $v_2$ and $v_3$, $v_2$ is adjacent to $v_4, \dots, v_k$, $v_3$ is adjacent to $v_4, \dots, v_k$, and there are no other adjacencies.
%
Let $A = \{v_1, \dots, v_k\}$ and $B = \{u_1,u_k, v_1 ,v_4, \dots, v_k\}$. Clearly, $A$ is an irregular dominating set of $G$. Let $w_1 = v_3$, $w_2 = v_2$, $w_3 = u_{k-1}, w_4 = u_{k-2}, \dots, w_k = u_2$. Thus, $V(G) \backslash B = \{w_1, \dots, w_k\}$. Note that $|N_{\bar{G}}(w_i) \cap B| = i$ for each $i \in [k]$. Thus, $B$ is an irregular dominating set of $\bar{G}$. Therefore, we have $\gamma_{ir}(G) \geq |A| = k$ and $\gamma_{ir}(\bar{G}) \geq |B| = k$, and hence the lower bounds in (i) and (ii) are attained.
 
Suppose that $n=6$, $u_1, u_2, u_3, v_1, v_2, v_3$ are the vertices of $G$, and $\{u_1,v_1\}, \{u_2,v_2\}$, $\{u_2,v_3\},\{u_3,v_1\},\{u_3,v_2\},\{u_3,v_3\}$ are the edges of $G$. Clearly, $\{v_1, v_2, v_3\}$ is an irregular dominating set of $G$, and $\{u_1, v_1, v_3\}$ is an irregular dominating set of $\bar{G}$. Thus, the lower bounds in (i) and (ii) are attained. 

Finally, suppose that $n=4$ and $G$ is the path $P_4 = ([4], \{\{1,2\}, \{2,3\}, \{3,4\}\})$. Then $\{1,3\}$ is an irregular dominating set of $G$, and $\{1,2\}$ is an irregular dominating set of $\bar{G} = ([4], \{\{2,4\}, \{4,1\}, \{1,3\}\})$. Thus, the lower bounds in (i) and (ii) are attained. 
\end{proof}


\end{document}